%% file: main.tex
\documentclass[runningheads]{llncs}

\usepackage{graphicx}
\usepackage{amsmath,amsfonts}
\usepackage{xspace}
\usepackage{xcolor}
\usepackage[ruled,linesnumbered]{algorithm2e}
  
  \SetCommentSty{commentstyle}
\usepackage[labelformat=simple]{subfig}
	
\usepackage{multirow}
\usepackage{hyperref}
\input{macros}

\begin{document}

\title{Reachability of weakly nonlinear systems using Carleman linearization\thanks{The first author is partly supported by Agencia Nacional de Investigaci\'{o}n e Innovaci\'{o}n (ANII), Uruguay.}}

\author{%
Marcelo Forets\inst{1}\orcidID{0000-0002-9831-7801}
\and \\
Christian Schilling\inst{2}\orcidID{0000-0003-3658-1065}
}
\authorrunning{M.\ Forets \and C.\ Schilling}

\institute{%
DMA, CURE, Universidad de la Rep\'ublica, Uruguay
\\
\email{mforets@gmail.com} \\
\and
University of Konstanz, Konstanz, Germany
\email{christian.schilling@uni-konstanz.de}
}

\maketitle

\begin{abstract}
In this article we introduce a solution method for a special class of nonlinear initial-value problems using set-based propagation techniques. The novelty of the approach is that we employ a particular embedding (Carleman linearization) to leverage recent advances of high-dimensional reachability solvers for linear ordinary differential equations based on the support function. Using a global error bound for the Carleman linearization abstraction, we are able to describe the full set of behaviors of the system for sets of initial conditions and in dense time.

\end{abstract}


\section{Introduction} \label{sec:intro}

\input{introduction}

\section{Preliminaries} \label{sec:preliminaries}

\input{preliminaries}

\section{Carleman linearization} \label{sec:carleman}

\input{carleman}

\section{Set propagation} \label{sec:method}

\input{method}

\section{Reachability algorithm} \label{sec:estimates}

\input{estimation}

\section{Evaluation} \label{sec:evaluation}

\input{evaluation}

\section{Conclusions} \label{sec:conclusion}

\input{conclusions}

\section*{Acknowledgments}
\input{acks}

\bibliographystyle{splncs04}
\bibliography{biblio}

\end{document}

%% file: macros.tex
\newcommand{\R}{\mathbb{R}}
\newcommand{\N}{\mathbb{N}}

\newcommand{\set}[1]{\mathcal{#1}}
\newcommand{\X}{\set{X}}
\newcommand{\Y}{\set{Y}}
\newcommand{\B}{\set{B}}
\newcommand{\range}[1]{[{#1}]} 

\newcommand{\norm}[1]{\Vert #1 \Vert}

\newcommand{\Idn}{\ensuremath{\mathbb{I}_{n}\xspace}}

\newcommand{\dd}{\mathrm{d}}

\newcommand{\Rs}{\ensuremath{\mathcal{R}}\xspace}
\newcommand{\post}{\ensuremath{\mathit{post}}\xspace}
\newcommand{\linearize}{\ensuremath{\mathit{linearize}}\xspace}
\newcommand{\error}{\ensuremath{\mathit{error}}\xspace}

\newcommand{\asgn}{\ensuremath{\leftarrow}\xspace}

\newcommand{\proj}[1]{\ensuremath{\pi_{#1}}\xspace}
\newcommand{\eps}{\ensuremath{\varepsilon}\xspace}

%% file: introduction.tex
We consider the problem of solving a system of nonlinear ordinary differential equations (ODEs) for a set of initial states. This is better known as \emph{reachability analysis}.
While for linear systems there exist very efficient algorithms \cite{GirardLGM06,LeGuernicG10,althoff2016combining,BogomolovFFVPS18}, reachability analysis for nonlinear systems remains a challenging problem.

Traditional approaches \cite{althoff2020set} include those based on Taylor models \cite{ChenAS13}, simulation \cite{Donze10}, or hybridization \cite{LiBB20}.
In this paper we present a new approach to this problem by transforming the nonlinear system into an infinite-dimensional linear system, which we then truncate. This truncated model approximates the original system.

More specifically, our approach is based on Carleman linearization, which is an established method in mathematical nonlinear control but differs from the above-mentioned approaches. The Taylor-model approach truncates an infinite Taylor polynomial, while we truncate a linear system. Hybridization approaches linearize smalls regions in the state space, while we linearize the whole system.

To achieve good accuracy, the truncation results in a high-dimensional linear system. To solve such systems, we leverage efficient reachability solvers based on the support function that have recently been developed.

Our approach can be used to obtain approximate solutions very quickly, but in an unsound way. Alternatively, using an error estimate, one can obtain a sound overapproximation. Under certain conditions (essentially weak nonlinearity), the error estimate converges, resulting in a precise approximation.

%
%
%

\paragraph{Contributions.}
This paper makes the following contributions:

\begin{itemize}
	\item We revisit Carleman linearization and explain how it can be used as a fast but unsound way to propagate sets through a nonlinear dynamical system.
	
	\item We extend the approach to a sound and practical reachability algorithm for dissipative nonlinear dynamical systems.
	
	\item We evaluate the algorithm in two case studies and discuss its strengths.
\end{itemize}

\paragraph{Related work.}

The original idea by Carleman \cite{carleman1932application,kowalski1991nonlinear} did not receive much attention for several decades.
Steeb showed that, while the nonlinear system and its infinite-dimensional embedding share the same analytic solutions, the embedding may admit additional non-analytic solutions \cite{Steeb89}.
Carleman linearization has since been applied successfully in control theory \cite{germani2005filtering,collado2008modified,mozyrska2008carleman,rauh2009carleman} and physics and chemistry \cite{gaude2001solving,HashemianA15}.

Several works provide bounds on the approximation error of the truncated linearized system \cite{forets2017explicit,liu2021efficient}. In this paper we use the error bound derived in \cite{liu2021efficient}.

An approach that is related to ours transforms a nonlinear system into a linear or polynomial system via a ``change of bases,'' using polynomials instead of Kronecker powers, and derives conditions under which this transformation preserves invariants \cite{Sankaranarayanan16}.

\paragraph{Outline.}

The next section recalls the mathematical basis used in this paper.
Section~\ref{sec:carleman} introduces the classic Carleman linearization.
In Section~\ref{sec:method} we describe how to propagate sets using Carleman linearization.
In Section~\ref{sec:estimates} we extend this approach to a reachability algorithm for dissipative nonlinear dynamical systems.
We evaluate the algorithm in Section~\ref{sec:evaluation} and conclude in Section~\ref{sec:conclusion}.

%% file: preliminaries.tex
In this section we summarize the mathematical prerequisites to make this paper self-contained. For a detailed derivation of the Carleman linearization procedure we refer to \cite{forets2017explicit}.

\subsection{Vectors, norms, and sets}

Let $\N = \{ 1,2,\ldots \}$ be the set of positive integers and $\R$ the set of real numbers, and for any $N \in \N$ we let $\range{N} := \{1, 2, \ldots, N\}$. $n$-dimensional vectors $x \in \R^n$ are understood as column vectors with components $x_i \in \R$, $i \in \range{n}$. Transposition is written $x^T$.
For any $x \in \mathbb{R}^n$ and $p \in \R_{\geq 1} \cup \{\infty\}$, $\norm{x}_p$ denotes the vector $p$-norm of $x$, with notable special cases $p=2$ (Euclidean norm) and $p=\infty$ (supremum norm). If $A = (a_{ij}) \in \R^{m \times n}$ is a matrix, $\norm{A}_p$ denotes the matrix norm induced by the vector $p$-norm, with notable special cases $p=2$ (spectral norm: the largest singular value) and $p=\infty$ (supremum norm: the maximum absolute row sum). We may abbreviate $\norm{\cdot}$ for $\norm{\cdot}_\infty$. See \cite{horn2012matrix} for precise definitions of these concepts.

For a compact, i.e. bounded and closed set $\X \subseteq \R^n$, $\norm{\X}_p$ denotes the maximum of $\norm{x}_p$ over all $x \in \X$. If $\X$ is polytopic (i.e. admits a representation as the finite intersection of half-spaces), its norm can be computed by a finite number of vector $p$-norm evaluations. Indeed, the map $x \to \norm{x}_p$ is convex and the maximum of a convex function over a polytope is attained at one of its vertices. However, computing the vertex representation of a polytope initially given by its half-space representation can be computationally expensive in dimensions higher than two (see \cite{kaibel2003}). A simpler rule applies if $\X$ is hyperrectangular (i.e., can be represented as an axis-aligned box with center $c \in \R^n$ and radius vector $r \in \R^n$). Then $\norm{\X}_p = \norm{c + D r}_p$ where $D = (D_{ij}) \in \mathbb{R}^{n \times n}$ is diagonal with matrix elements $D_{ii} = 1$ if $c_i \geq 0$ and $D_{ii} = -1$ otherwise, $i \in \range{N}$. We write $\B_r^n$ for the $n$-dimensional infinity-norm ball with radius $r$ centered in the origin. The projection of a set $\X$ to the first $k$ dimensions is denoted by $\proj{1:k}(\X)$.

\subsection{Support function}

A standard approach to operate with compact and convex sets in $\R^n$ is to use the \emph{support function} \cite{LeGuernic09}. The support function of $\X \subseteq \R^n$ along direction $d \in \R^n$, $\rho(d, \X)$, is the maximum of $d^T x$ over all $x \in \X$. In particular, if $\X$ is a polytope in half-space representation, its support function can be computed by solving a linear program (LP), and for certain classes of sets analytic formulas exist, which can be numerically evaluated in an efficient way. Such cases include hyperrectangular sets.
Since the support function distributes over Minkowski sums, i.e. $\rho(d, \X \oplus \Y) = \rho(d, \X) + \rho(d, \Y)$ for any pair of sets $\X$ and $\Y$, and since it holds that $\rho(d, M \X) = \rho(M^T d, \X)$ for any matrix $M \in \mathbb{R}^{n \times n}$, the support function has been successfully applied to solve linear set-based recurrences of the form $\X_{k+1} = M \X_k \oplus \Y_k$, either explicitly or implicitly by solving the recurrence only along a predefined number of directions \cite{LeGuernicG10,FrehseGDCRLRGDM11,BogomolovFFVPS18}. It is well-known that such recurrences are prevalent in reachability analysis of linear initial-value problems (IVPs), or nonlinear ones after some form of conservative linearization; see for example \cite{althoff2020set} and references therein.

\subsection{Kronecker product}\label{sec:kronecker}

For any pair of vectors $x \in \R^n$, $y\in \R^m$, their \textit{Kronecker product} is $w = x\otimes y = (x_1 y_1, \ldots, x_1 y_m, x_2 y_1,\ldots, x_n y_m)^T$, and the dimension is $\dim(w) = mn$. This product is not commutative. For matrices the definition is analogous: if $A \in \R^{m\times n}$ and $B \in \R^{p \times q}$, then $A \otimes B \in \R^{mp \times nq}$ and
\begin{equation*}
A \otimes B := \begin{pmatrix}
a_{11}B &\ldots&a_{1n}B \\
\vdots & & \vdots \\
a_{m1}B & \ldots &a_{mn}B 
\end{pmatrix}.
\end{equation*} 

The \textit{Kronecker power} $x^{\otimes i}$ of $x \in \mathbb{R}^n$ is a convenient notation to express all possible products of elements of a vector up to a given order:
\begin{equation*}
x^{\otimes i} :=  \underset{i \text{ times}}{\underbrace{x\otimes \cdots \otimes x}},\qquad x \in \R^n.
\end{equation*}
Note that $\dim(x^{\otimes i}) = n^i$, and each component  of $x^{\otimes i}$ is of the form $x_1^{\omega_1} x_2^{\omega_2} \cdots x_n^{\omega_n}$ for some multi-index $\omega \in \N^n$, $\norm{\omega}_1 = i$. For example, if $n=2$, the first two Kronecker powers are $x^{\otimes 1}= x = (x_1, x_2)^T$ and $x^{\otimes 2} = x \otimes x = (x_1^2,x_1 x_2,x_2 x_1,x_2^2)^T$. 
Further properties of Kronecker products
can be found in \cite{zhang2011matrix} and \cite{steeb2011matrix}.

%% file: carleman.tex
In this section we recall the classic Carleman linearization approach \cite{carleman1932application,kowalski1991nonlinear}.


Polynomial differential equations are an important class of nonlinear systems $x'(t) = f(x(t))$, $f : \R^n \to \R^n$, such that the coordinate functions $f_i : \R^n \to \R$ are multivariate polynomials. Many systems can be rewritten as polynomial vector fields by introducing
auxiliary variables, and any polynomial system is formally equivalent to a second-order system, possibly in higher dimensions -- for a proof of this statement and an algorithm to compute such transformation see \cite{forets2017explicit}. We can thus focus on quadratic ODEs without loss of generality. Consider the IVP for an $n$-dimensional quadratic ODE,
\begin{equation}\label{eq:quadraticODE}
    \frac{\dd x(t)}{\dd t} = F_1 x + F_2 x^{\otimes 2},
\end{equation}
with initial condition $x(0) \in \mathbb{R}^n$. Each $x_i(t)$, $i \in \range{n}$, is a function of $t$ over the interval $[0, T]$ where $T$ is the \emph{time horizon}. We assume that the matrices $F_1 \in \R^{n \times n}$ and $F_2 \in \R^{n\times n^2}$ are independent of $t$. Intuitively, $F_1$ (resp.\ $F_2$) is associated with the linear (resp.\ nonlinear) behavior of the dynamical system; thus $\norm{F_2}_2 / \norm{F_1}_2$ being small corresponds to \emph{weak nonlinearity} -- a concept we will use in a later section.

The \emph{Carleman linearization} (or \emph{Carleman embedding}) procedure begins by introducing a sequence of auxiliary variables $\hat{y}_j := x^{\otimes j}$, $j \in \N$. Differentiating such variables with respect to time, and repeatedly substituting~\eqref{eq:quadraticODE} into the derivatives of each $\hat{y}_j$ gives a formal equivalence with an infinite-dimensional linear system of ODEs \cite{kowalski1991nonlinear}. Truncation to order $N$ leads to a \emph{finite} linear IVP in the \emph{lifted} variables $\hat{y} := (\hat{y}_1, \hat{y}_2, \ldots, \hat{y}_N)^T$, namely
\begin{equation}\label{eq:LODE}
  \frac{\dd{\hat y}}{\dd{t}} = A \hat y, \qquad
  \hat y(0) = \hat y_{0},
\end{equation}
with initial condition $\hat y_{0} = (x_{0}, x_{0}^{\otimes 2}, \ldots, x_{0}^{\otimes N})^T$ and coefficients matrix $A$, which has the bi-diagonal block structure

\begin{equation}\label{eq:UODE}
A =
  \begin{pmatrix}
    A_1^1 & A_2^1 & 0 & 0 & \cdots & 0  \\
   0  & A_2^2 & A_3^2 & 0  & \cdots & 0 \\
    0 & 0 & A_3^3 & A_4^3 & 0 &  \vdots\\
    \vdots & \vdots & \vdots & \ddots & \ddots & 0 \\
    0 & 0 & \cdots & 0 & A_{N-1}^{N-1} & A_N^{N-1} \\
    0 & 0 & \cdots & 0 & 0 & A_N^N \\
  \end{pmatrix},
\end{equation}
where the $A^i_i$ and $A^i_{i+1}$, which we call \textit{transfer matrices}, have dimensions $n^i \times n^i$ and $n^i \times n^{i+1}$ respectively, and are defined by the formula 
\begin{equation*}
A^i_{i+i^\prime-1} = \sum_{\nu=1}^i   \overset{\text{i factors}}{\overbrace{\Idn \otimes \cdots \otimes \underset{\underset{\nu\text{-th position}}{\uparrow}}{F_{i^\prime}} \otimes  \cdots  \otimes \Idn }}.
\end{equation*}
for all $i \in \range{N}$ and where $i^\prime$ is either $1$ ($A^i_i$ is placed on the main diagonal) or $2$ ($A^i_{i+1}$ is placed on the upper diagonal) and where $\Idn$ is the identity matrix of order $n$. Note also that $A^1_1 = F_1$ and $A^1_2 = F_2$.
The dimension of \eqref{eq:LODE} is $n+n^2+\cdots+n^N=\frac{n^{N+1}-n}{n-1}=\mathcal{O}(n^N)$.



\paragraph{Running example.}



We illustrate the concepts described above in the simplest possible scenario.
Consider the logistic equation (a special case of~\eqref{eq:quadraticODE} for $n = 1$)
\begin{equation}\label{eq:logistic}
\dfrac{\dd x(t)}{\dd t} = rx(1 - x/K).
\end{equation}

This equation and related generalizations arise naturally in the context of population dynamics,
where $r > 0$ controls the initial rate of exponential growth, and $K > 0$ is the asymptotic equilibrium (the other equilibrium being $x = 0$).
We transform \eqref{eq:logistic} into the canonical scalar form \eqref{eq:quadraticODE}, namely $x'(t) = ax(t) + bx^2(t)$, via $a = r$ and $b = -r/K$.
Defining the auxiliary variables $\hat{y}_j = x^j$, $j \in \N$, we see that their first-order derivatives satisfy $\hat{y}_1' = x' = a\hat{y}_1 + b\hat{y}_2$, $\hat{y}_2' = 2x' x = 2a\hat{y}_2 + 2b\hat{y}_3$, etc. Hence the nonlinear ODE \eqref{eq:logistic} is equivalent to the (infinite) linear ODE
\begin{equation*}
    \hat{y}'_j = ja\hat{y}_j + jb\hat{y}_{j+1},\qquad j \in \N.
\end{equation*}

If we now fix the truncation order $N$, say, to $N = 4$, we obtain
\begin{equation*}
\dfrac{\dd \hat{y}(t)}{\dd t} = \begin{pmatrix}
    a & b & 0 & 0 \\
    0 & 2a & 2b & 0 \\
    0 & 0 & 3a & 3b \\
    0 & 0 & 0 & 4a \\
    \end{pmatrix} \hat{y},\qquad \hat{y}(0) = \begin{pmatrix} x_0 \\ x_0^2 \\ x_0^3 \\ x_0^4\end{pmatrix}.
\end{equation*}

\begin{figure}[tb]
	\centering
	\subfloat[Initial condition $x_0 = 0.5$.]{
		\centering
		\includegraphics[width=.48\textwidth,height=5cm,keepaspectratio]{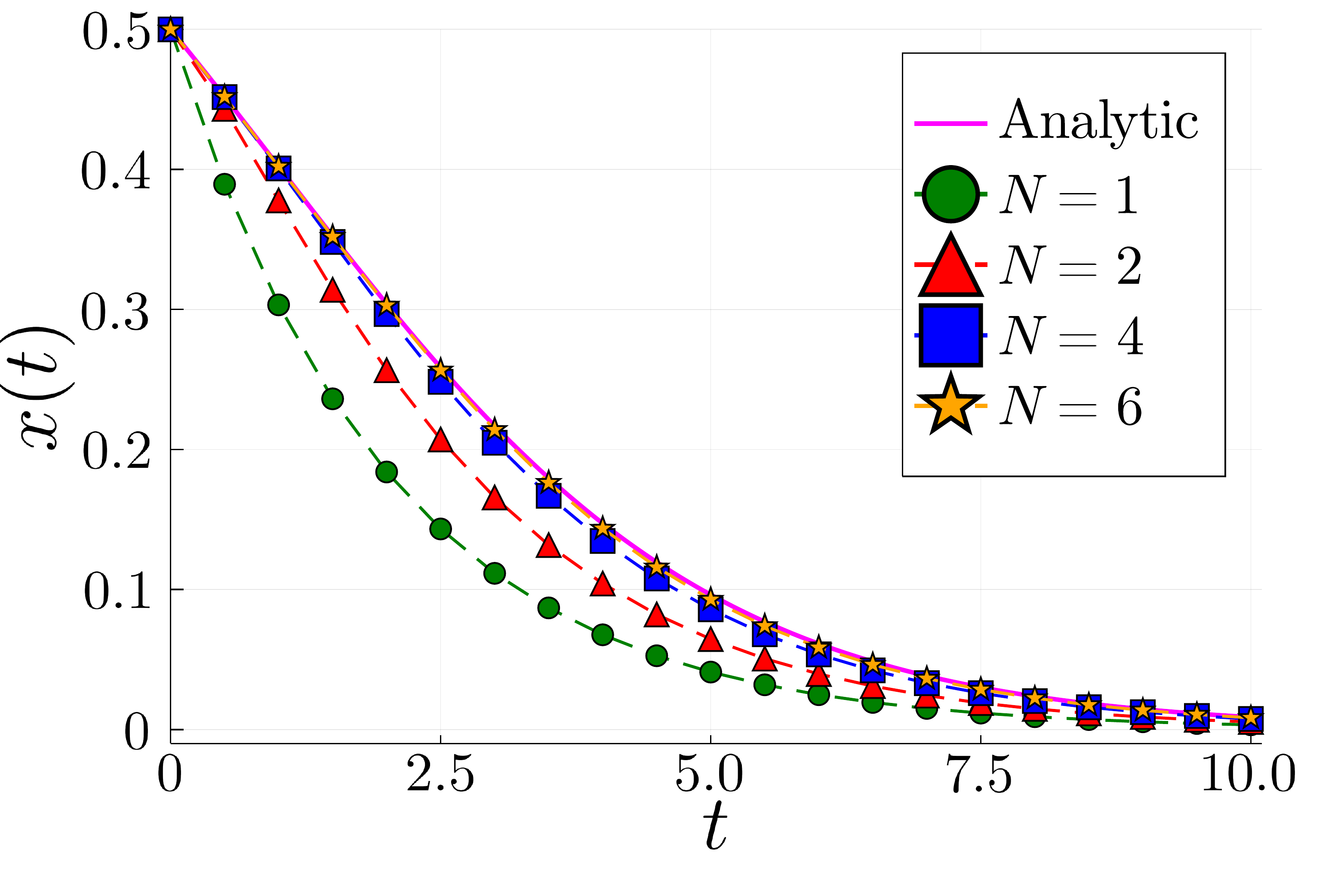} 
		\label{fig:logistic_nobloat_point}
	}
	\hspace*{-2mm}
	\subfloat[Initial condition $\X_0 = {[0.47, 0.53]}.$]{
		\centering
		\includegraphics[width=.48\textwidth,height=5cm,keepaspectratio]{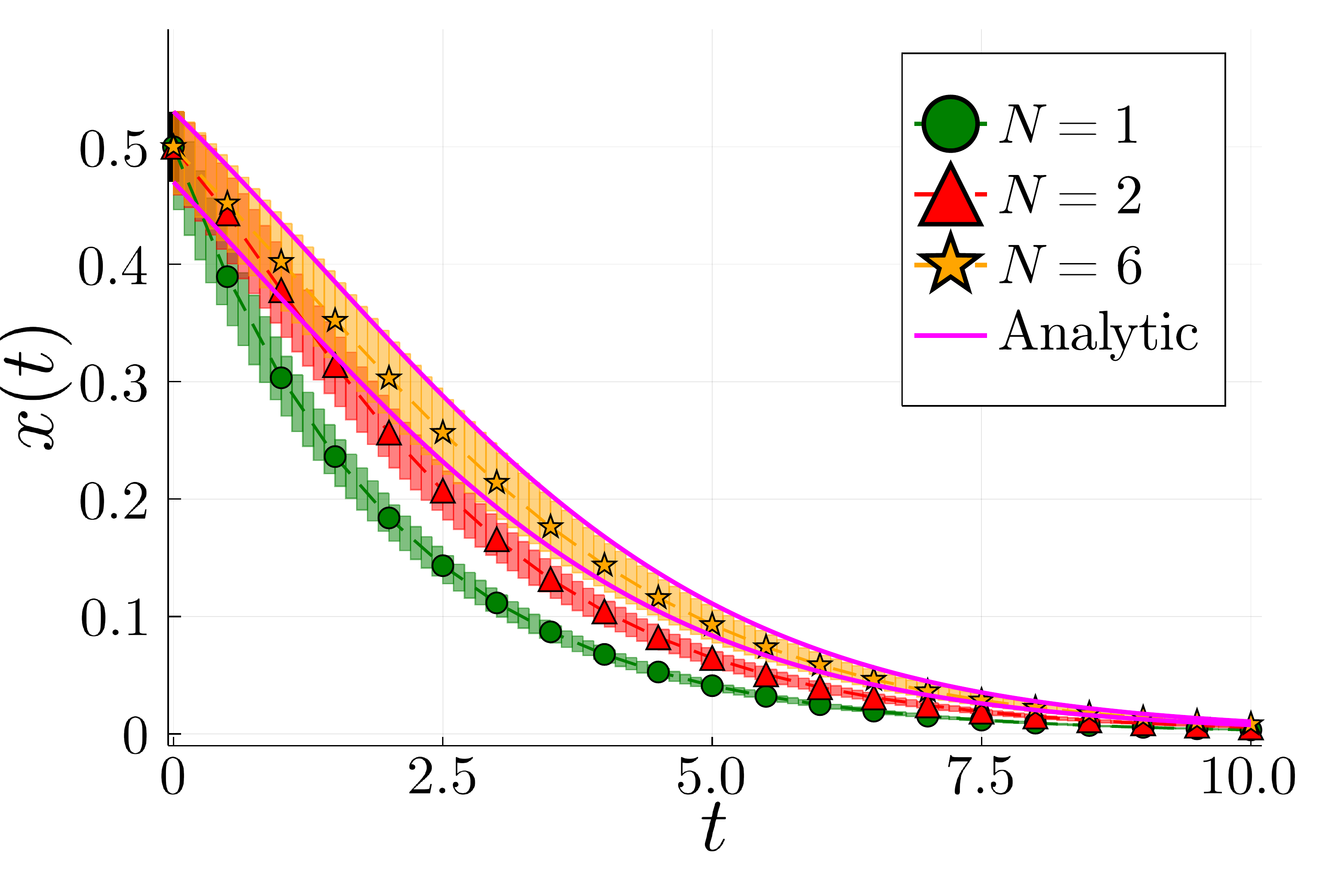}
		\label{fig:logistic_nobloat_set}
	}
	\caption{Solution of the original and of the truncated linearized IVPs for~\eqref{eq:logistic}.}
	\label{fig:logistic_nobloat}
\end{figure}

To estimate the quality of the approximation by finite truncation, we plot the solutions of \eqref{eq:logistic} from $x_0 = 0.5$ over a time horizon of $10$ and also plot the solution of \eqref{eq:LODE} for several choices of $N$ in Figure~\ref{fig:logistic_nobloat_point}. The model parameters are $r = -0.5$ and $K = 0.8$.
As can be seen, for increasing $N$ the solutions converge to the analytic solution, which in this case is known and given as
\[
	x(t) = \frac{x_0 a e^{a t}}{a + b (1 - e^{a t}) x_0}.
\]

Solving Eq.~\eqref{eq:LODE} requires computing the matrix exponential acting on the initial states $\hat{y}(t = k\delta) = e^{A k \delta} \hat{y}(0)$ at all times, which may be expensive for higher-dimensional systems. In the next section we introduce a method to propagate sets of initial conditions in dense time making use of the particular structure of the matrix \eqref{eq:UODE} using support function techniques. Theoretical estimates of the truncation error are considered in Section~\ref{sec:estimates}.


%
%

%% file: method.tex
In the previous section we saw how to transform a nonlinear IVP into an approximate linear IVP by Carleman linearization and truncation at a chosen order $N$. In this section we describe how this approach generalizes to IVPs whose initial condition is a set of states $\X_0 \subseteq \R^n$ described by a hyperrectangle. This is a common case, and hyperrectangular approximations can be computed efficiently.
We need to discuss two steps: how to transform $\X_0$ to the linear system and how to propagate sets of states for a linear IVP.

\smallskip

For the transformation of $\X_0$ we generalize the Kronecker product to sets with $\mathcal{X}^{\otimes i} := \{x^{\otimes i} \mid x \in \X\}$. For a hyperrectangle $\X$ we approximate $\mathcal{X}^{\otimes i}$ by applying the rules of interval arithmetic to each dimension \cite{moore2009introduction}. We note that one needs to carefully arrange the variables in order to obtain a tight solution. The arrangement consists of grouping the same variables of each monomial; for example, $x_1^2 x_2 x_1$ is evaluated using interval arithmetic as $x_1^3 x_2$ to avoid the dependency problem. 
To illustrate, consider the extension of the example in Section~\ref{sec:kronecker} to the hyperrectangle $\mathcal{X} = [0.9, 1.1] \times [-0.1, 0.1]$. Then $\mathcal{X}^{\otimes 1} = \mathcal{X} \subseteq \mathbb{R}^2$ and $\mathcal{X}^{\otimes 2} \subseteq [0.81, 1.21] \times [-0.11, 0.11] \times [-0.11, 0.11] \times [0.0, 0.01] \subseteq \mathbb{R}^4$.

\smallskip

There exist many algorithms to propagate a set through an IVP in a conservative way, i.e., the result overapproximates the true solution, in particular for linear IVPs \cite{Girard05,GirardLGM06,HanK06,LeGuernicG10,kaynama2011complexity,althoff2016combining,bak2017simulation,BogomolovFFVPS18,althoff2020set}. Most of these approaches first discretize the continuous-time system, for which the error can be made arbitrarily small by choosing a small discretization step $\delta$, and then propagate the sets in discrete time, which in certain cases can be done in an error-free way. We refer the reader to the above works for details about the discretization. Below we explain the second step because it is relevant for the later discussion.

Given a discretized linear IVP with discretized matrix $\Phi = e^{A \delta}$ and discretized initial condition $\hat{\X}_0$,
\begin{equation*}\label{eq:ivp_lin_disc}
	x_{k+1} = \Phi x_k, \quad x_0 \in \hat{\X}_0,
\end{equation*}
the set of reachable states is described by the \emph{flowpipe} $\bigcup_{k \geq 0} \Rs_k$ where the $\Rs_k := \Phi^k \hat{\X}_0$ is the \emph{reach set} for the time span $[k \delta, (k+1)\delta]$. In other words, a flowpipe is a sequence of reach sets given by the matrix powers of $\Phi$ applied to $\hat{\X}_0$. This computation scales to systems with hundreds of dimensions.

\paragraph{Example (cont'd).}

Consider again the logistic system. In Figure~\ref{fig:logistic_nobloat_set} we plot the flowpipes obtained for the different truncated approximations with an initial condition $\X_0 = [0.47, 0.53]$.

%% file: estimation.tex
In this section we discuss an error estimation that allows us to obtain a sound overapproximation of the states reachable by the original nonlinear system.

\subsection{Error bound}



We have yet to determine how the solutions of the truncated linear IVP~\eqref{eq:LODE} are related to those of the original nonlinear IVP~\eqref{eq:quadraticODE}. To formulate this relation precisely, we introduce some notation. The error of the $j$-th block of variables is defined as $\eta_j(t) := x^{\otimes j}(t) - \hat{y}_j(t)$, which is the difference between the Kronecker power of the solution of~\eqref{eq:quadraticODE} and the projection of the solution of~\eqref{eq:LODE} onto the corresponding block of variables of the lifted $\hat{y}$. We are mostly interested in the first block, i.e., $j=1$, since $x(t) = x^{\otimes 1}(t)$, and the truncation error corresponds to upper bounding the quantity $\norm{\eta_1(t)} \leq \eps(t)$ for some error function $\eps(t)$ to be determined. Ideally, for fixed $t$ the error function should decrease sufficiently fast for increasing order $N$, so we can use low orders in practice, typically 2 to 6.
In \cite{forets2017explicit} the authors derived explicit error bounds for the linearization, i.e., a function $\eps(t)$ that only depends on the initial condition and the norms of the matrices $F_1$ and $F_2$. However, that approach is too conservative since $\eps(t)$ diverges in finite time -- even in cases when the solution of the linearized system~\eqref{eq:LODE} is converging.

Crucial to the present article, the authors in \cite{liu2021efficient} discovered that, by imposing an assumption on the class of quadratic problems considered, an arbitrary-time and exponentially convergent error formula holds. There are two main assumptions: 1)~linear terms dominate over nonlinear ones (\emph{weak nonlinearity}) and 2)~nonlinear effects play a prominent role during a finite time span, after which only linear terms matter (\emph{linear dissipation}). These definitions are formalized below. In the following we assume that the eigenvalues of $F_1$ in~\eqref{eq:quadraticODE} are sorted (counting multiplicities) such that $\Re(\lambda_n) \leq \cdots \leq \Re(\lambda_1)$, where $\Re(\lambda)$ is the real part of $\lambda$.

\begin{definition} \label{hy:weak}
System~\eqref{eq:quadraticODE} is said to be \emph{weakly nonlinear} if the ratio
\begin{equation}\label{eq:defR}
    R := \dfrac{\norm{x_0}  \norm{F_2}}{\vert \Re(\lambda_1) \vert}
\end{equation}
satisfies $R < 1$.
\end{definition}

\begin{definition}\label{hy:dissipative}
System~\eqref{eq:quadraticODE} is said to be \emph{dissipative} if $\Re(\lambda_1) < 0$ (i.e., the real part of all eigenvalues is negative).
\end{definition}

The conditions $\Re(\lambda_1) < 0$ and $R < 1$ ensure arbitrary-time convergence.

\begin{theorem}[{\cite[Corollary 1]{liu2021efficient}}]\label{thm:error}
Assuming that~\eqref{eq:quadraticODE} is weakly nonlinear and dissipative, the error bound associated with the linearized problem \eqref{eq:LODE} truncated at order $N$ satisfies
\begin{equation}
\norm{\eta_1(t)} \leq \eps(t) := \norm{x_0} R^N (1 - e^{\Re(\lambda_1) t})^N, \label{eq:error}
\end{equation}
with $R$ as defined in~\eqref{eq:defR}. This error bound holds for all $t \geq 0$.
\end{theorem}

\subsection{Obtaining a sound set-propagation algorithm}

The interesting aspect of~\eqref{eq:error} is that we can enclose all possible behaviors of a nonlinear problem for a hyperrectangular initial condition $\X_0 \subseteq \R^n$ in two steps: first, propagating the solutions of the high-dimensional linear system \eqref{eq:LODE} forward in time using a suitable linear reachability technique; in a second step, enlarging the solution (a sequence of reach sets $\Rs_j$ with associated time span $\Delta t = [t, t + \delta]$ for some $\delta > 0$) by taking the Minkowski sum with a ball of radius $r := \max(|a|, |b|)$ where $[a, b]$ is the interval-arithmetic evaluation of $\eps(\Delta t)$. Moreover, the truncation error converges to zero for increasing $N$ and, as we will see in  the experiments, typical values of $N$ do not have to be prohibitively large to obtain reasonable approximation bounds.

\begin{theorem}
Given a flowpipe, consider any $n$-dimensional reach set $\Rs_j$, $j \geq 0$, and its associated time span $\Delta t = [t, t + \delta]$. Let $\overline{\Rs}_k$ be the true set of reachable states in the time span $\Delta t$ and $r$ as defined above. Then we have $\overline{\Rs}_k \subseteq \Rs_k \oplus \B_{r}^n$.
\end{theorem}

%



This allows us to present a sound reachability method as shown in Algorithm~\ref{algo:reach}.
Crucially, we see that in Line~\ref{line:bloating} we only require the reach sets $\Rs_j$ in the first $n$ dimensions. Thus it suffices to compute these sets in a ``sparse'' way in Line~\ref{line:linreach}. We can use an algorithm based on the support function for $\post$ to achieve that. In our implementation we use the algorithm from \cite{LeGuernicG10}, which takes as input a set of direction vectors in which the reach sets are evaluated.

\begin{algorithm}[t]
    \caption{Reachability algorithm}
    \label{algo:reach}
    \KwIn{%
    $\X_0 \subseteq \R^n$: hyperrectangular initial states;
    $F_1, F_2$: system matrices;
    \goodbreak
    $N$: truncation order;
    $T$: time horizon;
    \post: algorithm to compute a flowpipe for linear systems \\
    }
    \vspace*{1mm}
    \KwOut{flowpipe overapproximating the reachable states until $T$}
    \vspace*{2mm}
    %
    $A, \hat{\X}_0 \asgn \linearize(\X_0, F_1, F_2, N)$ \tcp*{Carleman linearization}
    %
    $(\Rs_0, \dots, \Rs_m) \asgn \post(y' = Ay, y(0) \in \hat{\X}_0)$ \tcp*{flowpipe for linear system} \label{line:linreach}
    \For{$j \asgn 0$ \KwTo $m$}{%
      $\eps \asgn \error(\Rs_j, \X_0, F_1, F_2, N)$ \tcp*{linearization error}
      $\Rs_j^\eps \asgn \proj{1:n}(\Rs_j) \oplus \B_\eps^n$ \tcp*{enlarged reach set} \label{line:bloating}
    }
    \Return{$(\Rs_0^\eps, \dots, \Rs_m^\eps)$}
\end{algorithm}


\subsection{Reevaluation of the error term}

For dissipative systems, while the solution of the linear system may converge to zero, the corrected term including the error estimate may not. This observation leads to the idea of reevaluating the error estimate after some time $t^*$, since for fixed $F_1$ and $F_2$, a decreasing $\norm{x_0}$ leads to a smaller value $R$ which, in turn, reduces the error estimate $\eps(t)$. This is, however, nontrivial because by the time one reevaluates, the past error estimate must be taken into account and thus the new state estimate at $t^*$ may already be too pessimistic. In the evaluation we apply such a reevaluation manually.


%% file: evaluation.tex
In this section we study two models that have also been used in \cite{liu2021efficient}, but we repeat them here to make this article self-contained.
In the first model we evaluate all aspects outlined in the present article including the error bounds. In the second model we demonstrate that even if the assumptions for the error bounds do not apply, we can still obtain solutions of useful accuracy.

For comparison we compute an overapproximation of the reachable states for the original nonlinear systems using a Taylor-model (TM) approach implemented in JuliaReach~\cite{BogomolovFFPS19,TaylorModels.jl,TaylorIntegration.jl,TaylorSeriesJOSS}, with the default parameters (Taylor polynomials with spatial and temporal expansions of orders $2$ and $8$ respectively), which generally has high precision.
To evaluate the flowpipe for the linear system in Algorithm \ref{algo:reach}, we use the $2n$ directions $\pm e_i$ for $i \in \range{n}$, where $e_i$ is the unit vector in dimension $i$, which corresponds to the outer hyperrectangular approximation of $\proj{1:n}(\Rs_j)$. Note that the number of directions, $2n$, is independent of the truncation order $N$.
Interval-arithmetic computations are performed using the Julia library IntervalArithmetic.jl \cite{IntervalArithmetic.jl}, and for set-based computations we use LazySets.jl \cite{LazySets.jl}. The code and scripts to run these problems is available online. \footnote{\href{https://github.com/JuliaReach/RP21_RE}{github.com/JuliaReach/RP21\_RE}}

\subsection{Epidemic model (SEIR)}

There exist several widely used models of population dynamics that generalize the logistic model from Section~\ref{sec:carleman} \cite{brauer2012mathematical}.
We consider the popular \emph{SEIR} epidemic model with data on the early spread of the COVID-19 disease from \cite{seir20}. A population $P$ is divided into four compartments: susceptible ($P_S$), exposed ($P_E$), infectious ($P_I$), and recovered ($P_R$). An individual is initially susceptible and becomes exposed/infected with rate $r_\text{tra}$. The latent time before an exposed individual becomes infectious themselves is $T_\text{lat}$. Finally, an infectious individual recovers after time $T_\text{inf}$. New individuals are added to the population with rate $\Lambda$. We also consider a vaccination with rate $r_\text{vac}$ \cite{ZamanKJ08}. The system of ODEs is:

\begin{align*}
	\frac{\dd P_S}{\dd t} &= -\Lambda \frac{P_S}{P} - r_\text{vac} P_S - r_\text{tra} P_S \frac{P_I}{P} + \Lambda
	\\
	\frac{\dd P_E}{\dd t} &= -\Lambda \frac{P_E}{P} - \frac{P_E}{T_\text{lat}} + r_\text{tra} P_S \frac{P_I}{P} \\
	\frac{\dd P_I}{\dd t} &= -\Lambda \frac{P_I}{P} + \frac{P_E}{T_\text{lat}} - \frac{P_I}{T_\text{inf}} \\
	\frac{\dd P_R}{\dd t} &= -\Lambda \frac{P_R}{P} + r_\text{vac} P_S + \frac{P_I}{T_\text{inf}}
\end{align*}

In this model we assume that $P = P_S + P_E + P_I + P_R$ remains constant, so we need not model $P_R$. The corresponding $F_i$ matrices thus simplify to

\begin{align*}
	F_1 = \begin{pmatrix}
		-\frac{\Lambda}{P} - r_\text{vac} & 0 & 0 \\
		0 & -\frac{\Lambda}{P} - \frac{1}{T_\text{lat}} & 0 \\
		0 & \frac{1}{T_\text{lat}} & -\frac{\Lambda}{P} - \frac{1}{T_\text{inf}}
	\end{pmatrix},
	F_2 = \begin{pmatrix}
		0 &\ 0 &\ -\frac{r_\text{tra}}{P} &\ 0 &\ 0 &\ 0 &\ 0 &\ 0 &\ 0 \\
		0 &\ 0 &\ \frac{r_\text{tra}}{P} &\ 0 &\ 0 &\ 0 &\ 0 &\ 0 &\ 0 \\
		0 &\ 0 &\ 0 &\ 0 &\ 0 &\ 0 &\ 0 &\ 0 &\ 0
	\end{pmatrix}.
\end{align*}

Since $F_1$ is triangular, $\Re(\lambda_1) = -\frac{\Lambda}{P} - \min\{r_\text{vac}, \frac{1}{T_\text{lat}}, \frac{1}{T_\text{inf}}\}$. We also have
$\norm{F_2} = \frac{\sqrt{2} r_\text{tra}}{P}$ and
$\norm{\X_0} \leq P$.
Thus we can estimate

\[
	R = \dfrac{\norm{\X_0} \frac{\sqrt{2} r_\text{tra}}{P}}{\frac{\Lambda}{P} + \min\{r_\text{vac}, \frac{1}{T_\text{lat}}, \frac{1}{T_\text{inf}}\}} \leq \dfrac{\sqrt{2} r_\text{tra}}{\frac{\Lambda}{P} + \min\{r_\text{vac}, \frac{1}{T_\text{lat}}, \frac{1}{T_\text{inf}}\}}.
\]


The time scale is measured in days. We use the same parameters as in \cite{liu2021efficient}: a population of $P = 10^7$, $\Lambda$ is small (here: $\Lambda = 1$), hence the constant term is disregarded in the analysis, $T_\text{lat} = 5.2$, $T_\text{inf} = 2.3$, $r_\text{tra} = 0.13$ days$^{-1}$, and $r_\text{vac} = 0.19$ days$^{-1}$. We choose $\X_0 = [6e6, 3e5, 3.7e6] \oplus \B_{1e5}^3$, which results in $R \approx 0.68$ and $\Re(\lambda_1) \approx -0.19$ and thus Theorem~\ref{thm:error} is applicable.



\begin{figure}[tb]
	\centering
	\subfloat[Reach sets $\Rs_j$ without error estimate.]{
		\centering
		\includegraphics[width=.48\textwidth,height=5cm,keepaspectratio]{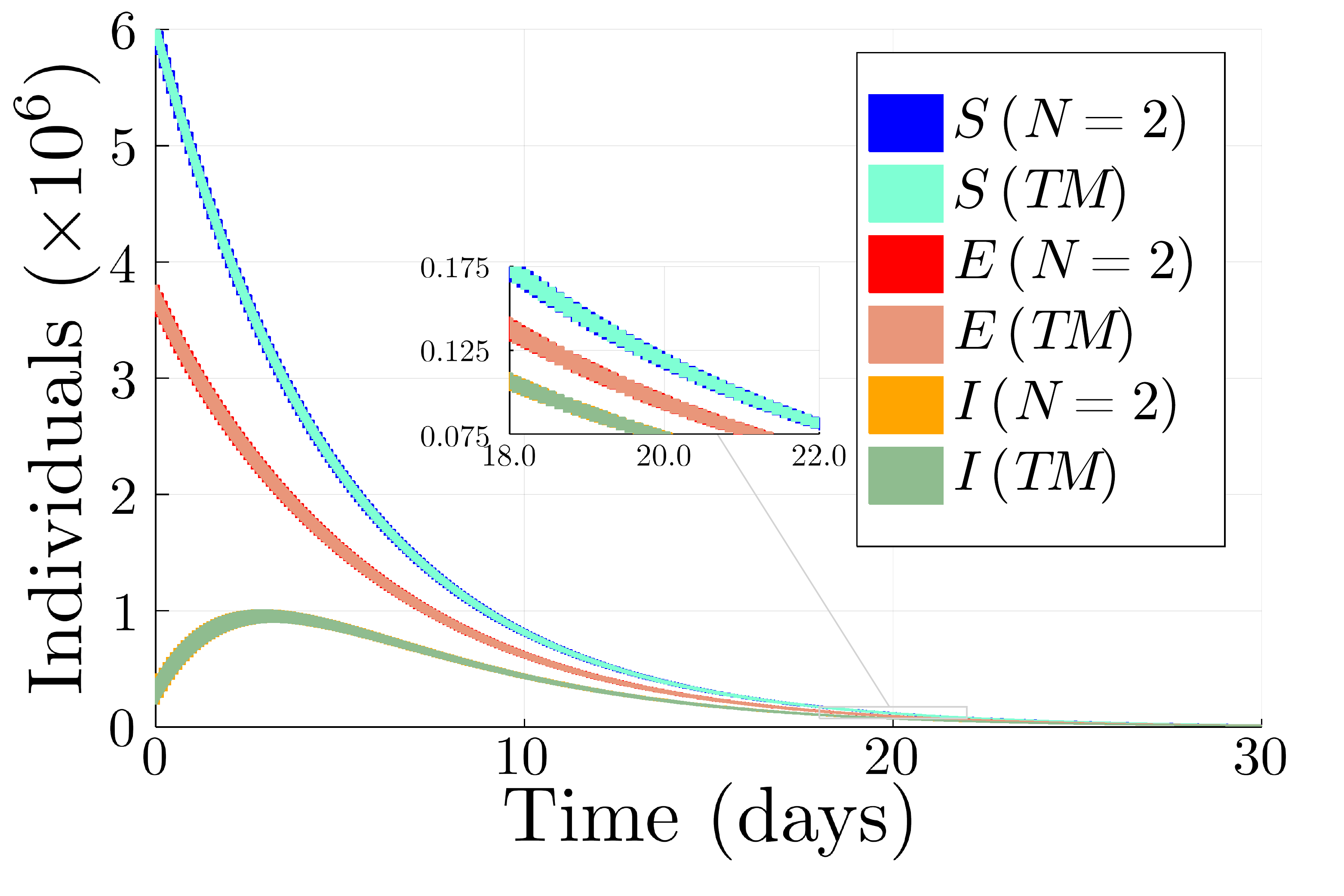}
		\label{fig:seir_noerror}
	}
	\hspace*{-2mm}
	\subfloat[Reach sets $\Rs_j^\eps$ including error estimate with a reevaluation at $t = 4$ (cf.\ zoom).]{
		\centering
		\includegraphics[width=.48\textwidth,height=5cm,keepaspectratio]{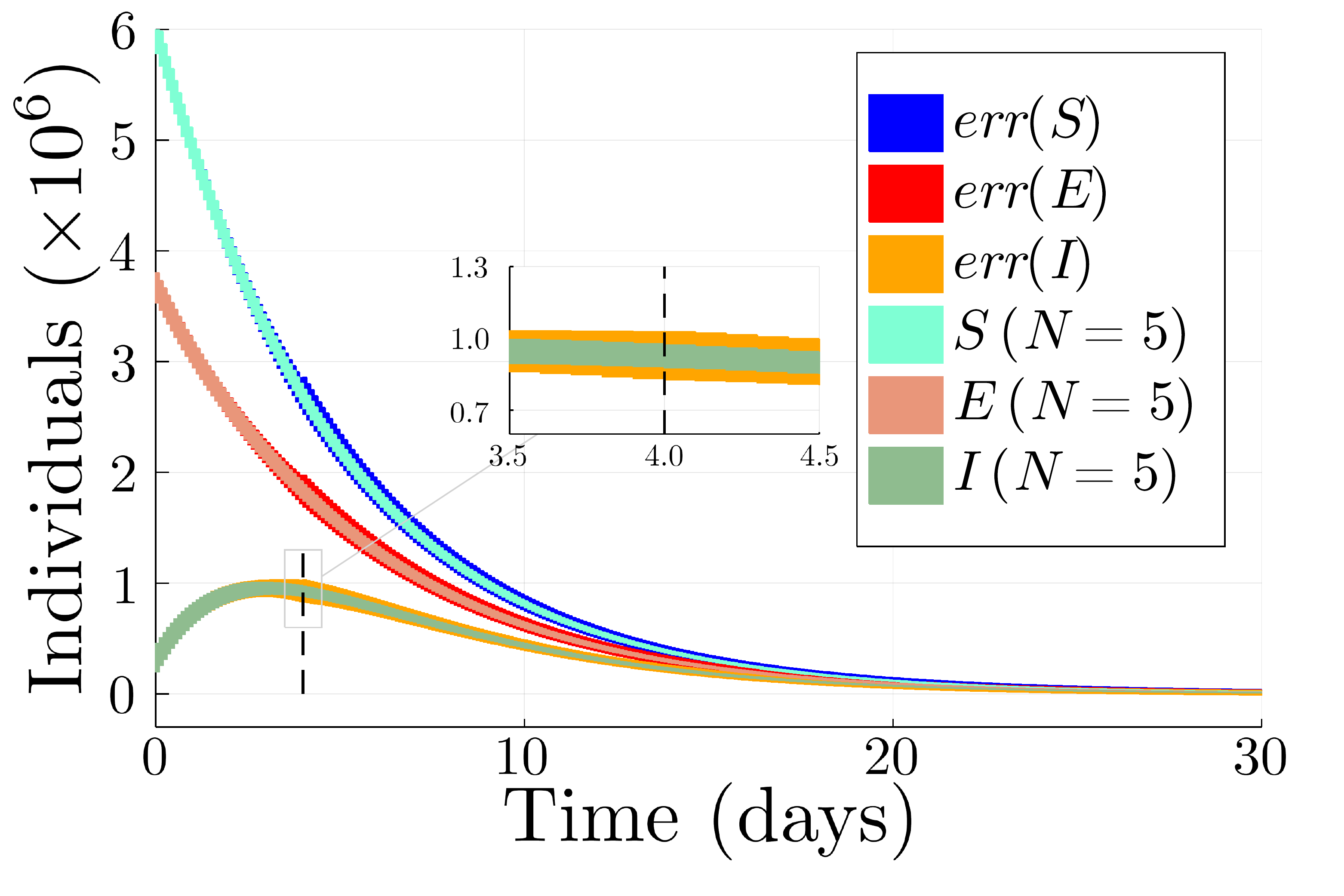}
		\label{fig:seir_error}
	}
	\caption{Results for the SEIR model.}
	\label{fig:seir}
\end{figure}

\begin{table}[tb]
	\centering
	\begin{tabular}{c @{\hspace*{3mm}} l @{\hspace*{3mm}} l @{\hspace*{3mm}} c @{\hspace*{3mm}} l @{\hspace*{3mm}} l}
		& \multicolumn{2}{c}{SEIR model} && \multicolumn{2}{c}{Burgers model} \\
		\cline{2-3} \cline{5-6} \\
		& no error bound & incl.\ error bound && initial point & initial set \\
		\hline
		TM & \multicolumn{2}{c}{$6.14$\,s} && \multicolumn{1}{c}{$0.88$\,s} & \multicolumn{1}{c}{$0.91$\,s} \\
		\hline
		\multirow{2}{*}{Carleman} & $N = 2$: $0.006$\,s & $N = 5$: $0.185$\,s && $N = 2$: $0.0065$\,s & $N = 2$: $0.0067$\,s \\
		&  &  && $N = 3$: $0.24$\,s & $N = 3$: $0.29$\,s
	\end{tabular}
	\caption{Run times for the SEIR model and the Burgers model obtained for the Taylor-model (TM) approach and the Carleman linearization with different truncation orders $N$.}
	\label{tab:runtimes}
\end{table}

The analysis results without and with conservative error estimate are plotted in Figure~\ref{fig:seir}, where we used the discretization step $\delta = 0.1$. We can see that the non-conservative Carleman approximation is precise even for the small value $N = 2$. However, the error estimate is too conservative for such small value of $N$ thus it is not plotted. However, using $N=5$ the error estimate improves significantly, but only until around time $t = 4$; this is due to the large values in $\X_0$. At $t = 4$ we reevaluate the estimate. Since the state and with it the norm has changed, the new error estimate is more optimistic and converges quickly. The run times are given in Table~\ref{tab:runtimes}.


%
%
%
%

%
%
%
%

\subsection{Burgers partial differential equation}

\begin{figure}[tb]
	\centering
	\subfloat[Point initial condition ($w = 0$).]{
		\centering
		\includegraphics[width=.482\textwidth,height=5cm,keepaspectratio]{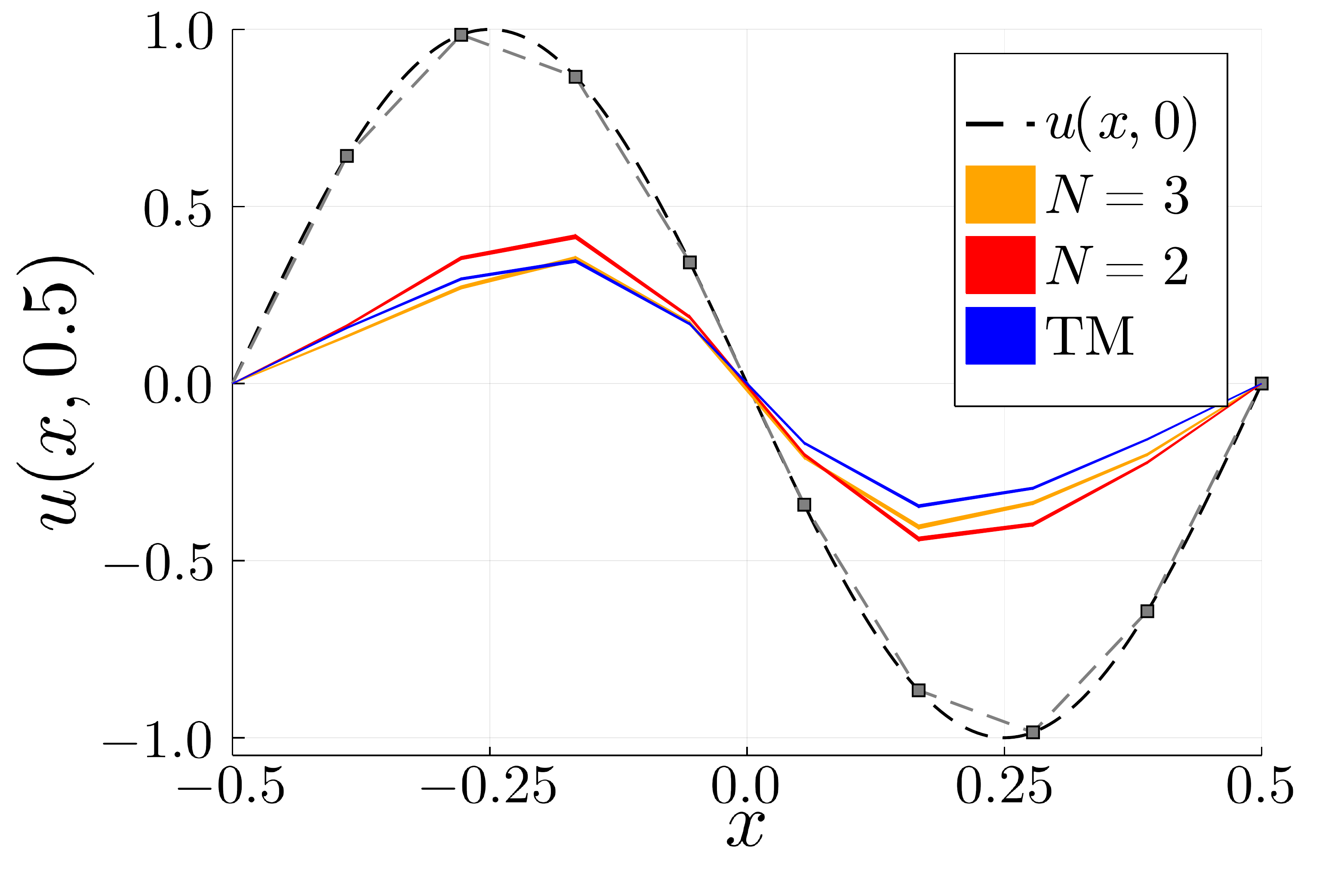}
		\label{fig:burgers_point}
	}
	\hspace*{-2mm}
	\subfloat[Distributed initial condition ($w = 0.06$).]{
		\centering
		\includegraphics[width=.482\textwidth,height=5cm,keepaspectratio]{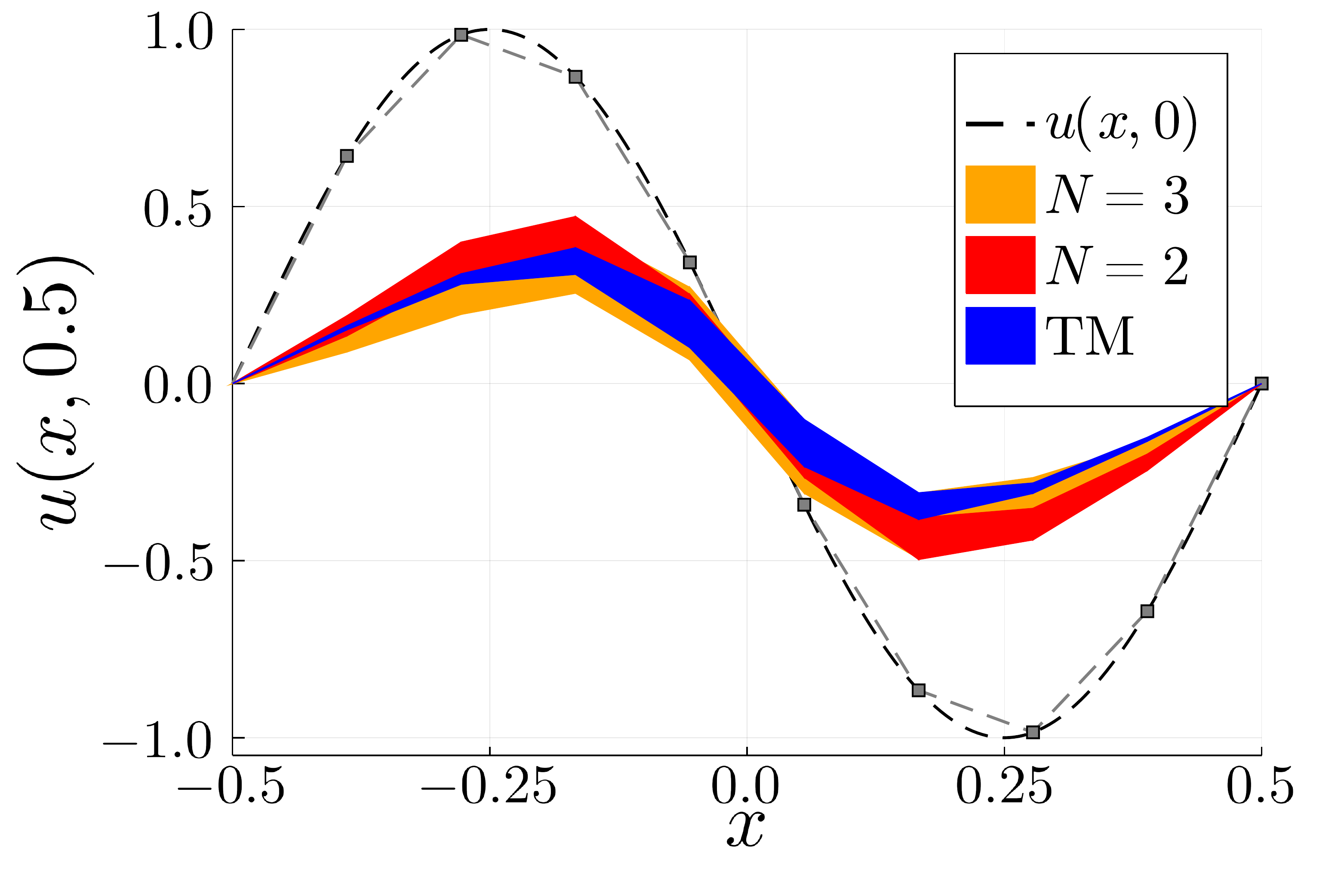}
		\label{fig:burgers_set}
	}
	\caption{Results for the Burgers model at $t = 0.5$ with initial condition width $w$.}
	\label{fig:burgers}
\end{figure}

We study a model arising from the discretization of a partial differential equation (PDE). Consider the viscous Burgers equation to model convective flow \cite{Burgers48}
\begin{equation*}
	\partial_t u + x \partial_x u = \nu \partial_x^2 u.
\end{equation*}

We use the following model parameters: viscosity $\nu = 0.05$, domain length $L_0 = 1$, and $U_0 = 1$. We consider the initial condition $u(x, 0) = -U_0 \sin(2 \pi x / L_0)$ on the domain $x \in \pm L_0/2$ and Dirichlet conditions $u(x, 0) = 0$ at the boundaries. We distribute this initial condition to a set by keeping the end points fixed and enlarging the initial point to some width $w = 0.06$. For the PDE discretization we use central differences obtaining the coupled differential equations
\begin{equation}\label{eq:burgers_discretized}
	\partial_t u_i = \nu \frac{u_{i+1} - 2 u_i + u_{i-1}}{\Delta x^2} - \frac{u_{i+1}^2 - u_{i-1}^2}{4 \Delta x}.
\end{equation}

We use $n_x = 10$ points and $\Delta x = L_0 / (n_x - 1)$. Eq.~\eqref{eq:burgers_discretized} has the form of \eqref{eq:quadraticODE} that we need to apply Carleman linearization. We obtain $\Re(\lambda_1) \approx -0.488 < 0$ but $R \approx 18.58$, i.e., $R$ as defined in Eq.~\ref{eq:defR} is not smaller than one. Although the theoretical error bounds from Theorem~\ref{thm:error} are not applicable here, it is interesting to observe that the set-based solution is reasonably accurate with respect to the solution obtained for the original nonlinear system.
In Figure~\ref{fig:burgers} we plot the results at $t = 0.5$. For the linear reachability algorithm we used the step size $\delta = 0.01$. We can see that we still obtain good approximations that decrease exponentially by incrementing the truncation $N$. The run times are given in Table~\ref{tab:runtimes}.

%% file: conclusions.tex
In this paper we have presented a reachability method that abstracts nonlinear terms into a higher-dimensional space such that the evolution is approximately linear. The main advantage of the method is that we can leverage recent set propagation techniques that are specialized to high-dimensional linear ODEs. However, the method does not apply to general nonlinear systems but requires weak nonlinearity, i.e., the relative norm of the nonlinear term should be smaller than that of the linear term. Under such limitations, the presented method outperforms other reachability methods because linear reachability in high dimension can be solved efficiently.

This work can be extended in several ways.
First, we can consider time-dependent terms; an error bound is derived in \cite{liu2021efficient}.
Second, in our experimental evaluation we observed that manually reevaluating the error bound can improve the precision if the norm of the states shrinks (which should generally happen for dissipative systems). It would be interesting to automate this process.
Third, the reachability analysis can be accelerated, e.g., using Krylov methods to work more efficiently in high dimensions. A more challenging direction is to devise a new reachability algorithm that exploits the structure of the linearized system.

%% file: acks.tex
The first author acknowledges fruitful discussions with Amaury Pouly and Goran Frehse. We are grateful to the anonymous reviewers for helpful comments.